\documentclass{IEEEtran4PSCC}
\usepackage{xcolor}
\ifCLASSINFOpdf
   \usepackage[pdftex]{graphicx}
   \usepackage{multicol}
\else
   \usepackage[dvips]{graphicx}
\fi
   \usepackage[cmex10]{amsmath}
   \usepackage{multicol}
   \usepackage{amssymb}
   
\makeatletter
\let\old@ps@headings\ps@headings
\let\old@ps@IEEEtitlepagestyle\ps@IEEEtitlepagestyle
\def\psccfooter#1{%
    \def\ps@headings{%
        \old@ps@headings%
        \def\@oddfoot{\strut\hfill#1\hfill\strut}%
        \def\@evenfoot{\strut\hfill#1\hfill\strut}%
    }%
    \def\ps@IEEEtitlepagestyle{%
        \old@ps@IEEEtitlepagestyle%
        \def\@oddfoot{\strut\hfill#1\hfill\strut}%
        \def\@evenfoot{\strut\hfill#1\hfill\strut}%
    }%
    \ps@headings%
}
\makeatother

\begin{document}
\title{Modeling of Network Constraints in Large-scale Capacity Expansion Optimization of Power Grids}

\author{
\IEEEauthorblockN{
Yu Weng\IEEEauthorrefmark{1}, 
Lara Booth\IEEEauthorrefmark{2}, 
Priya L. Donti\IEEEauthorrefmark{2}, 
Ruaridh Macdonald\IEEEauthorrefmark{1}}
\IEEEauthorblockA{\IEEEauthorrefmark{1} MIT Energy Initiative, Cambridge, USA}
\IEEEauthorblockA{\IEEEauthorrefmark{2} Department of Electrical Engineering \& Computer Science, and Laboratory for Information \& Decision Systems, \\
Massachusetts Institute of Technology, Cambridge, USA\\ \{wengyu22, lcbooth, donti, rmacd\}@mit.edu}
}


\maketitle

\begin{abstract}
Capacity expansion modeling plays a critical role in optimizing the deployment of new generation, storage, and transmission, typically at national and regional levels. To support long-term planning, these models consider a large set of energy technologies and policies, along with decades of weather and demand data. Realistic capacity expansion models thus become high-dimensional optimization problems, with hundreds of millions of variables and constraints, which are challenging to solve. A common strategy to address this complexity is to omit non-linear, non-convex AC optimal power flow (ACOPF) constraints and instead use linearized power balance equations or transport formulations. While these simplifications improve tractability, they limit our understanding of how power flow and the physical properties of power networks impact investment decisions across generation, storage, and transmission infrastructure. This paper addresses this gap by extending the GenX capacity expansion model to incorporate fixed point theorem-based network constraints. These embed ACOPF-based considerations while maintaining the tractability of the planning model, nearly preserving the dimensionality of the transport formulation and incurring only modest runtime increases.
This approach is much cheaper than embedding ACOPF directly, making it appropriate for large-scale capacity planning problems. We compare our approach to the original transport-based GenX model as well as a non-linear, non-convex version that incorporates the full ACOPF constraints, for a case study of the ISO New England grid.

\thanksto{\noindent Corresponding author: Ruaridh Macdonald (rmacd@mit.edu).\\
This publication was based upon work supported by the U.S. Department of Energy’s Office
of Energy Efficiency and Renewable Energy (EERE) under the Hydrogen Fuel Cell Technology Office, Award Number DE-EE0010724
}

\end{abstract}


\section{Introduction}

Power systems are in urgent need of capacity expansion and network updates to enable the siting and integration of low-carbon energy resources while also meeting rising demand from data centers and electrified loads \cite{eia_steo_2026_03, iea_electricity_2026}. Capacity expansion models (CEMs) play a critical role in optimizing the long-term planning behind energy infrastructure deployments, typically at national or state levels. 
To enable holistic planning, CEMs must be able to optimize over a wide range of energy technologies and diverse operational, public policy, and planning constraints.
Recent research has shown that CEMs must consider decades of operational data \cite{ruggles2024planning, bhatt2025missing} at high temporal \cite{levin2024high} and spatial resolutions \cite{qiu2024decarbonized, krishnan2016evaluating-594, serpe2025importance} to produce robust electricity grid designs. \textcolor{black}{Models must also include} numerous interdependent technologies, resulting in high-dimensional optimization problems \textcolor{black}{that are generally not tractable. To allow these problems to be solved, CEMs use reduced temporal and/or spatial resolutions and apply simplifications} such as constraint linearization, resource aggregation, and decomposition techniques. Although these strategies reduce computational burden, they often come at the cost of model fidelity. Comparative studies have evaluated these trade-offs over a range of applications \cite{van2022impacts}.

One of the most significant simplifications in CEMs is the omission of AC optimal power flow (ACOPF) constraints. Nearly all large-scale planning models omit these constraints \cite{wogrin2020assessing} due to the computational burden imposed by their non-linearity and non-convexity. 
While there are some exceptions, specifically in the transmission expansion literature \cite{garcia2025ac, torres2014expansion, khanpourempirical}, most CEMs incorporating generation and transmission expansion solely include a simple linearized power balance. A few studies have examined how this simplification affects investment decisions in generation and transmission infrastructure. Among them, \cite{neumann2022assessments, lee2025canopi} examined the impact of approximated power flow and transmission losses via various DC-OPF formulations, and \cite{van2022impacts} developed a CEM using a second-order cone programming approximation of ACOPF. \cite{recht2024considering} incorporated reactive power into a CEM. All found that more advanced representations of transmission produced different planning results compared to simpler models. However, none directly compared their results to a CEM using full ACOPF constraints, leaving open the question of exactly how much these intermediate approximations improve model accuracy.

This work aims to advance our understanding of how AC power flow and physical properties of power networks impact investment decisions, and proposes a scalable methodology for incorporating higher-fidelity physics. Our contributions are:

\begin{enumerate}
    \item We integrate AC power flow constraints into GenX \cite{genxmainref, Bonaldo_GenX_2025}, a large-scale capacity modeling framework. This enables more physically faithful assessment of how AC grid constraints influence investment decisions. 
    \item In a case study of New England, we show that incorporating ACOPF constraints \textcolor{black}{leads to additional generation capacity being installed and a more even distribution of generation across nodes, due to the explicit consideration of reactive power demand, generators having reduced active power dispatch capacity, and networks having reduced active power transfer capacity and greater losses.}
    \item While our ACOPF-CEM results show the importance of considering AC physics, including them significantly increases CEM runtime. We thus propose an approach that incorporates AC power flow-related considerations using fixed point (FX) theorem-based network constraints. These FX constraints can readily be embedded into large-scale linear CEMs, reproducing ACOPF-like behavior while maintaining computational tractability. We validate our approach against the non-linear, non-convex version of GenX with full AC power flow constraints.
    \item Finally, we show that the computed FX network constraints are robust to out-of-sample conditions, such as altered demand shapes and fuel price variations, and thus suitable for planning problems under uncertainty.
  
\end{enumerate}
\section{Modeling of Capacity Expansion Planning}\label{sec_cems}

 GenX is an open source CEM intended for studying and optimizing electricity grids. A full description of GenX and its capacity expansion formulation is available in \cite{genxmainref, Bonaldo_GenX_2025}. In compact form, the original mixed-integer linear program (MILP) version of GenX can be described as follows:
 
\begin{align}
\min_{x, y} \quad &  \sum_z^\mathcal{Z} \sum_t^T \sum_g^\mathcal{G} (C_\mathcal{I} y+ C_\mathcal{O}x  ) \label{eq_GenXModel} \\
        \text{subject to} \quad & A_\mathcal{O}x + B_\mathcal{I}y \leq r , \label{GenXcnst_IO} \\ &F_\mathcal{I}y \leq d,  \label{GenXcnst_I} \\ & g_\mathcal{O}(x,y) \leq 0,   \label{GenXcnst_NL}\\ & x, y \geq 0, x,y \in R  \nonumber \label{GenXcnst_xy}.
\end{align}

\noindent The goal is to minimize the total system cost resulting from planning decisions $y$, and dispatch decisions $x$, while satisfying exogenous electricity demands and policy constraints. The cost parameter $C$ includes two components: investment costs, denoted by subscript $\mathcal{I}$, which account for expenditures such as new generation capacity and transmission lines; and operational costs, denoted by $\mathcal{O}$, which capture expenses like fuel consumption and penalties. The sets $\mathcal{T}$, $\mathcal{Z}$, and $\mathcal{G}$ represent the time periods, zones, and technologies considered in the system, respectively, with $t$, $z$, and $g$ referring to specific elements within those sets. Constraints \eqref{GenXcnst_IO} contain all operational constraints, including capacity limits of all resources, policy-related constraints, and technology-specific constraints. Examples of the latter are ramping limits and shut-on/shut-off constraints for thermal generators, or state of charge and charging/discharging limits for energy storage. Constraints \eqref{GenXcnst_I} relate to the investment and retirement decisions of all generation, storage, and transmission assets. Integrality constraints are included in \eqref{GenXcnst_NL}, which enable GenX to be configured as a mixed-integer linear optimization model to handle binary and integer variables related to unit commitment and reserves. 

GenX makes use of a simplified formulation of transmission operations in order to reduce the runtime of the underlying MILP and allow for large-scale models to be solved. The set of variables $x$ includes a signed power flow variable on each transmission line, $l$, with flow direction between zones determined by a network map. The line flow variables are constrained to be at most the available transmission capacity of that line \eqref{GenXcnst_IO}. Active power line losses are set to be a fixed fraction of the magnitude of the flows on each line. We refer to this transmission formulation as the \textit{transport formulation}.

\section{Full AC power flow and feasibility constraints in CEMs}
\label{sec_ACOPF}

This section introduces our approach to building AC power flow into existing capacity expansion models. To capture the physical behavior of the power grid, we integrate AC optimal power flow constraints--using the formulation given in Power Grid Lib \cite{pglib}, adapted for multiple time steps--into the existing GenX optimization model. These constraints govern the feasibility of generator dispatches, power transfers, and voltage levels across zones and time periods. For all $t \in \mathcal{T}$, the AC power flow constraints are as follows:
\begin{align}
\angle V_{\mathrm{ref},t} = 0 \quad & \\
\underline{\Delta \theta}_l \leq \angle V_{i,t} - \angle V_{j,t} \leq \overline{\Delta \theta}_l \quad &\forall l \in \mathcal{L} \\
\underline v_z \leq |V_{z,t}| \leq \overline v_z \quad &\forall z \in \mathcal{Z} \\
\underline S_{g} \leq S_{g,t} \leq \overline S_{g} \quad &\forall g \in \mathcal{G} \\
S_{l,t} = \left( Y^*_l - i \tfrac{B^c_l}{2} \right) \tfrac{|V_{i,t}|^2}{|T_l|^2} - Y^*_l \tfrac{V_{i,t} V^*_{j,t}}{T_l} \quad &\forall l \in \mathcal{L} \\
S^\mathrm{rev}_{l,t} = \left( Y^*_l - i \tfrac{B^c_l}{2} \right) |V_{j,t}|^2 - Y^*_l \tfrac{V^*_{i,t} V_{j,t}}{T^*_l} \quad &\forall l \in \mathcal{L} \\
|S_{l,t}| \leq \overline S_l \quad \text{and} \quad |S^\mathrm{rev}_{l,t}| \leq \overline S_l \quad &\forall l \in \mathcal{L} \\
\sum_{g \in \mathcal{G}_z} S_{g,t} = \sum_{l \in \mathcal{L}_z} \left(S_{l,t} + S^\mathrm{rev}_{l,t} \right) \quad &\nonumber \\ + \; D_{z,t} + Y^{\mathrm{sh}}_z |V_{z,t}|^2 \quad &\forall z \in \mathcal{Z}. 
\end{align}

\noindent Our AC formulation uses the following notation:
\begin{itemize}
    \item $\mathcal{T}, \mathcal{Z}, \mathcal{G}, \mathcal{L}$: sets of time periods, zones, generators, and transmission lines, respectively.
    \item $V_{z,t} \in \mathbb{C}$: complex voltage at zone $z$ and time $t$.
    \item $S_{g,t} = P_{g,t} + jQ_{g,t}$: complex power generation of generator $g$ at time $t$.
    \item $S_{l,t}, S^\mathrm{rev}_{l,t} \in \mathbb{C}$: complex power flows on line $l$ at time $t$ in forward and reverse directions.
    \item $\underline{S}_g, \overline{S}_g \in \mathbb{R}_{\geq 0}$: apparent power limits of generator $g$.
    \item $\underline{v}_z, \overline{v}_z \in \mathbb{R}_{>0}$: voltage magnitude bounds at zone $z$.
    \item $\underline{\Delta\theta}_l, \overline{\Delta\theta}_l \in \mathbb{R}$: angle difference bounds over line $l$.
    \item $Y_l \in \mathbb{C}$: series admittance of line $l$.
    \item $B^c_l \in \mathbb{R}$: line charging susceptance of line $l$.
    \item $T_l \in \mathbb{C}$: complex tap ratio of line $l$; 1 if not a transformer.
    \item $D_{z,t} \in \mathbb{C}$: complex net demand at zone $z$ and time $t$.
    \item $Y^{\mathrm{sh}}_z \in \mathbb{C}$: shunt admittance at zone $z$.
\end{itemize}

We incorporate AC power flow into GenX by linking AC active generation variables with GenX’s generator dispatch variables and augmenting zonal power balance with the active power flows at the corresponding end of each incident line--forward flows for lines originating in the zone and reverse flows for lines terminating in the zone--as well as the shunt conductance term for each zone. Transmission losses are represented through the sum of the forward and reverse active power flows on each line. GenX's fixed fraction power loss from the transport formulation is removed to avoid double-counting losses. Apparent power flow on each line is constrained by the available transmission capacity in GenX. Note that, unlike in traditional ACOPF, we do not introduce a quadratic or piecewise linear dispatch cost term in the CEM objective, but instead maintain the original linear operational cost term. \textcolor{black}{No relaxation methods are applied for the ACOPF constraints. Our incorporation of reactive power into GenX is similar to previous attempts \cite{wogrin2020assessing, recht2024considering}, but these other studies did not use full ACOPF constraints. In our work, the resulting version of GenX is a non-convex, mixed-integer non-linear program.}

\begin{figure}[!t]
    \centering
\includegraphics[width=0.85\linewidth]{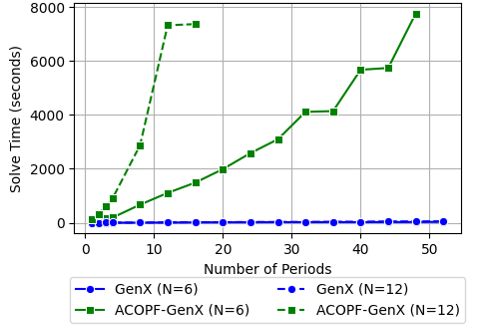}
    \vspace{-1em}
\caption{\textcolor{black}{Solve time comparison of a case study of the ISO New England grid, with and without ACOPF constraints. Each run has a different number of spatial zones ($N=6, 12$) and representative weeks (x-axis; 1,2,3, and multiples of 4 until 52). All runs had a 3 hour time limit.}}
    \label{fig:solvetime}
\end{figure}

\textcolor{black}{As we will show in Section \ref{sec_results},} directly integrating ACOPF 
enhances the accuracy and physical fidelity of the capacity expansion problem. \textcolor{black}{However, it also significantly increases runtime. Fig. \ref{fig:solvetime} compares the solve time of a case study of the ISO New England grid, with and without ACOPF. 
(Solver and hardware details are in Section~\ref{sec_results}.)
Thirty-two versions of the case were run: with either 6 or 12 zones (i.e., aggregated buses) and between 1 and 52 representative weeks. When using the transport formulation, we see that the runtime increases linearly and only slightly with problem size, from 10s to 100s for 1 to 52 representative periods, respectively, for the 12-zone case. By contrast, adding ACOPF causes GenX's runtime to increase significantly with the number of zones and exponentially with the number of representative weeks. All cases were run with a 3 hour time limit, which were exceeded by the 6-zone, 52-period and 12-zone, 20-period cases. 
Notably, these problems are lengthy to solve despite being relatively small. Recent studies have shown that investment models must have tens of zones \cite{krishnan2016evaluating-594, serpe2025importance} and consider decades of data \cite{ruggles2024planning, bhatt2025missing, levin2024high} to produce robust results. 
In other words, integrating ACOPF directly is likely to be computationally prohibitive for realistic-scale planning problems. While advances in non-linear solvers may improve the situation somewhat \cite{pacaud2024gpu}, overall, more tractable approaches are needed to integrate network physics into CEMs.}

\section{Network Feasibility Constraints in CEMs}

 To mitigate the computational burden of directly using ACOPF constraints, this section instead proposes fixed-point theorem-based (FX) network feasibility constraints that tractably represent AC-related considerations. 
 These FX constraints are integrated directly into the original MILP version of GenX that uses the transport transmission formulation. 
 Specifically, we maintain the same decision variables and constraints as the original transport GenX formulation (rather than adding additional variables as for ACOPF in Section~\ref{sec_ACOPF}). We then add additional FX constraints into the transport formulation, to limit decision variables to values that are consistent with AC feasibility. Our FX constraints contain a correction factor to account for the fact that CEMs use aggregated power network parameters (e.g., aggregating multiple buses and lines) and the fact that existing CEMs like GenX are written with active power--not apparent power--in mind. Together, this allows us to retain key physical characteristics of AC power systems while nonetheless maintaining computational tractability.

 In the rest of this section, we first introduce background on fixed-point certificates for AC power flow.
 We then derive our CEM-adapted certificate, accounting for network aggregation and active-power-only formulations.
 Finally, we describe the procedure for embedding the derived certificate within a CEM.
 
\subsection{Fixed-point theorem and derived certificate}
\label{subsec:fx-standard}

The proposed FX network feasibility constraints are built upon a feasibility set that preserves the physical limitations of power networks in planning problems. This feasibility set defines allowable power injections that satisfy both steady-state AC network stability and transmission line capacity limits \cite{weng2022asymmetrically,sampath2022voltage}. We compute this feasibility set using an approach based on the advanced fixed-point theorem-based certificate in \cite{weng2021fixed}. This section briefly summarizes the key ideas of the fixed-point theorem-based certificate and the feasibility set.

\subsubsection{Fixed-point theorem} 

In mathematics, fixed-point theorems specify conditions under which a mapping $F(\cdot)$ has at least one point $x$ satisfying $x = F(x)$ \cite{agarwal_meehan_oregan_2001}. Writing a system in the form $x=F(x)$ maps the variables into themselves, and provides a lens for assessing the solvability of equations. Pertinently, fixed-point approaches can be used to assess the existence of solutions to the AC power flow equations.

In contrast to the family of fixed-point methods proposed in \cite{nguyen2018constructing}, \emph{Kantorovich's fixed-point theorem} is closely related to the Newton-Raphson (NR) method and is therefore well-suited for integration with AC power flow equations. This theorem specifies three conditions on the initial point $x_0$ and the Jacobian $J(x)$ of the mapping $F(x)$: 
\begin{align}
    &\operatorname{det}(J(x_0)) \neq  0 \nonumber \\
    & \| (J(x_0))^{-1} (J(x)-J(y)) \| \leq L \| x-y \|,\\
    &\| (J(x_0))^{-1} (F(x_0)-b)  \| = \eta \leq \frac{1}{2L}.\nonumber 
\end{align}
When these conditions hold, the theorem guarantees a unique solution within the admissible region and ensures that NR initialized at $x_0$ converges to it. (Here, $b,L,\,\text{and}\;\eta$ are constants. Throughout this work, the operator $\| \cdot \| $  denotes the 1-norm.) 

However, applying Kantorovich’s theorem to power systems is not without obstacles. The standard power-flow relation $S = V I^*$ is not an analytical function, due to complex conjugation. As a result, the Jacobian $J(x)$ does not conform to the form required by the theorem and must be reformulated. The reformulation must also preserve the original physical properties of the power flow equations and ensure a non-conservative feasibility region for admissible injections. Solutions to these challenges were presented in \cite{weng2021fixed}.

\subsubsection{Fixed-point theorem-based certificate for power systems}

The certificate derived in \cite{weng2021fixed} extends Kantorovich's fixed-point theorem to the AC power flow equations to produce a sufficient condition for steady-state stability of power systems. For a power network operating at a nominal point with voltage $V_{nom}$ under power injections $S_{nom}$, satisfying $S = V I^*$, the condition can be written as:
\begin{align}\label{C:Kan-PL}
\| \zeta(S_{nom}) \| \|\zeta(S-S_{nom})\| 
&+\frac{\| J^{-1}_{nom}\Delta(S-S_{nom})\|}{2\| J^{-1}_{nom}  \|^2} \nonumber\\
&\leq \frac{1}{4 \| J^{-1}_{nom}  \|^2},
\end{align}
where 
\begin{equation} \label{eq_Zeta}
    \zeta(S)=[V_{nom}]^{-1} Z [V^*_{nom}]^{-1} [S^*],
\end{equation}
\begin{align} \label{eq:jnom}
 J_{nom}=\begin{bmatrix}
[\mathbf{1_n}] & \zeta(S_{nom}) \\
\zeta^*(S^*_{nom}) & [\mathbf{1_n}] 
\end{bmatrix},
\\
 \Delta(S)
=\begin{bmatrix} \label{eq:deltaS}
	[\zeta(S)\mathbf{1_n}]&  \zeta(S) \\
	{\zeta(S)} & [\zeta^*(S)\mathbf{1_n}]
  \end{bmatrix}.
\end{align}
Here, $Z$ is the system impedance matrix,
and the notation $[\cdot]$ denotes a diagonalization operator: given a vector, $[\cdot]$ forms the diagonal matrix with that vector on its diagonal. 
The certificate parameters in equations~\eqref{eq_Zeta}--\eqref{eq:deltaS} serve to bound feasible deviations from the fixed point. 
Under condition ~\eqref{C:Kan-PL}, a unique equilibrium is guaranteed and the network is in steady-state stability for any power injection $s$ satisfying either inequality \eqref{C:Kan-PL}.  Moreover, NR iterations initialized at $(V_{nom},S_{nom})$ converge to this unique solution. Note that this unique solution can only be obtained when $(V_{nom},S_{nom})$ is well-defined. If the system is already non-stable, the power injections and voltages at this non-stable point are not usable for the certificate above. In addition, the particular choice of nominal point may lead to condition~\eqref{C:Kan-PL} being more or less conservative, as well as affecting generalization to scenarios that depart significantly from the nominal point.
These equations also assume a static power grid topology and impedance matrix; future work will explore extending the certificate to accommodate varying topologies and impedances.

Conditions \eqref{C:Kan-PL} is straightforward to embed in optimization models at low computational cost, while the Jacobian and $Z$ matrix preserve the physical characteristics of AC power flow. No simplifying assumptions on the topology, parameters, or operating conditions are introduced. These properties provide the theoretical foundation for introducing network feasibility constraints within a CEM.
\subsection{Adapting FX network feasibility constraints for CEMs}
\label{sec:fxconstraints}


The features of CEMs, as compared to power flow and stability analyses, lead to practical differences in the optimization formulation and the computation of certificate parameters and linear constraints. The key differences of the formulation between CEMs and standard AC-based optimization include:

\begin{itemize}
    \item \textit{Limited spatial resolution and zonal numbers.} Manipulating the system impedance matrix is inexpensive, as the number of zones in CEMs is typically $<100$. 
    
    
    \item \textit{No fixed bus-type of nodes.} In CEMs, all zones are modeled as mixtures of generators and loads, with no differentiation between generator-bus vs. load-bus types. The ACOPF formulation in Section~\ref{sec_ACOPF} reflects this choice.
    
    \item \textit{Aggregated network parameters.} In CEMs, line parameters and transmission limits of the original power network are aggregated at the zonal level. Equivalent parameters must be calculated for the ACOPF constraints to be valid.

    \item \textit{Calculation units.} CEMs operate in engineering units rather than per-unit values. All variables must be rescaled to ensure the fixed-point formulation remains valid.
\end{itemize}
Given these details, we derive an alternative set of FX conditions to embed within the CEM as follows. 

For simplicity of analysis, we consider a nominal point with no power injections, \textit{i.e.}, $S_{nom} = 0$. Thus, there is no current flow and the voltage is $V_{nom}=1 \angle 0$ in p.u. From~\eqref{eq_Zeta}--\eqref{eq:jnom}:
\begin{equation*}
    \zeta(S_{nom}) = 0, \; J_{nom}=\begin{bmatrix}
\mathbf{1_n} & \mathbf{0_n} \\
\mathbf{0_n} & \mathbf{1_n} 	
\end{bmatrix}.
\end{equation*}
The certificate \eqref{C:Kan-PL} can then be simplified as:
\begin{equation} \label{C_Sis0}
    \| \Delta(S) \| \leq 1/2.
\end{equation}
This certificate defines admissible zonal  injections that remain consistent with AC network feasibility.

We note that the certificate in \eqref{C_Sis0} does not directly apply to an aggregated CEM, because the CEM uses network parameters and transmission limits that have been reduced to zonal quantities. 
In addition, the CEM uses only active power variables, not apparent power.
We thus modify the 
certificate~\eqref{C_Sis0} to account for aggregation and to express the constraints in terms of active rather than apparent power.
Specifically, we introduce a correction factor to account for the former, and use power factors to account for the latter.

Let $\mathbf{S}_0 \in \mathbb{C}^{Z \times T}$ denote the complex zonal power injections from a 
solution to CEM with ACOPF (Section~\ref{sec_ACOPF}),
and let $\mathbf{S}_0^l \in \mathbb{C}^{L \times T}$ denote the corresponding complex line flows. 
(Note that $S_{nom}$ and $\mathbf{S}_0$ may be different points.)
For each time step $t$, we define $\mathbf{s}_0(t) \triangleq \mathbf{S}_0(:,t) \in \mathbb{C}^{Z}$ and $\mathbf{s}_0^l(t) \triangleq \mathbf{S}_0^l(:,t) \in \mathbb{C}^{L}$. Define the apparent-power magnitude operator $\mathcal{M}_{ag}(S)=|S|=\sqrt{P^2+Q^2}$ and the power-factor operator $\varrho(S)=P/|S|$.
Using these quantities, we introduce a margin parameter, $\lambda$, that rescales the admissible region based on the reference ACOPF-CEM solution:
\begin{equation}
\lambda = \left\| \Delta \left( \mathcal{P}_r \mathbf{S}^l_{0, \max}(S_{\mathrm{cap}})^+ \right) \right\|, 
\qquad 
\mathcal{P}_r = \mathbf{S}_0 (\mathbf{S}_0^l)^+,
\label{eq_lambda}
\end{equation}
where $S_{\mathrm{cap}} \in \mathbb{R}^{L}$ is the vector of line apparent-power capacity limits, $\mathbf{S}^l_{0, \max}$ is the complex line flow at which the maximum apparent-power magnitude in $\mathbf{S}_0^l$ occurs, and $(\cdot)^+$ represents the pseudoinverse. Intuitively, $\lambda$ measures the admissible margin at the reference ACOPF-CEM solution relative to a benchmark upper flow limit. This factor is computed at the reference ACOPF-CEM solution and remains fixed in subsequent FX-constrained optimization.

Using this margin, we define the FX constraints as
\begin{align}\label{eq:Crr_and_bounds}
    \mathcal{C}_{rr} |\mathbf{S}|/|\mathbf{S}_0| \le \lambda,
\end{align}
where $\mathbf{S}$ is the candidate zonal apparent-power injection to be constrained. Here, $\mathcal{C}_{rr}$ does not define a second, independent margin. Instead, it determines how conservatively candidate operating points are constrained against the reference-point margin $\lambda$ in the aggregated CEM. We consider three variants:
\begin{align*}
\mathcal{C}_{rr}^{\mathrm{T}}(t) &= \left\| \Delta \Big(\mathcal{P}_r \mathbf{s}_0^{l}(t)/ \mathcal{M}_{ag}\big(\mathbf{s}_0^{l}(t)\big)\Big) \right\|, \\
\mathcal{C}_{rr}^{\mathrm{M}}(t) &= \left\| \Delta \big(\mathcal{P}_r \mathbf{s}_0^{l}(t)/ \overline{\mathcal{M}_{ag}}(\mathbf{S}_0^{l}) \big) \right\|, \\
\mathcal{C}_{rr}^{\mathrm{L}}(t) &= \left\| \Delta \big(\mathcal{P}_r \mathbf{s}_0^{l}(t)(S_{\mathrm{cap}})^+ \big) \right\|.
\end{align*}
These variants correspond to different ways of enforcing the reference-point margin. The tight bound uses the hourly line flows from the reference solution, the medium bound uses the maximum line flow over all time steps, and the loose bound uses the line capacities. This creates a trade-off between conservatism and flexibility: tighter bounds provide greater confidence in AC feasibility, while looser bounds admit solutions closer to the transport-model optimum. Note the loose bound does not reduce to the transport formulation, as the FX constraints still restrict deviations from the nominal point.

Conventional CEMs optimize active power, whereas the FX construction above is written in terms of apparent power. We therefore use the reference-solution power factor to express the constraints in active-power form. Let $\boldsymbol{\varrho}^{\mathbf{S}_0}$ denote the elementwise power factor associated with $\mathbf{S}_0$. Then the compact FX-CEM constraint~\eqref{eq:Crr_and_bounds} becomes
\begin{equation}
\mathcal{C}_{rr} \big(|\mathbf{P}|/ \boldsymbol{\varrho}^{\mathbf{S}_0} \big) / \mathcal{M}_{ag}(\mathbf{S}_0)\leq \lambda.
\label{eq_FXC}
\end{equation}

\noindent $\mathbf{P}$ denotes the active power constrained in the planning problem. In GenX, we replace $\mathbf{P}$ by $\mathcal{E}_{z,t}$, the net zonal active power injection and power scale factor $\varsigma$ (MW/GW), then the final version of the FX-CEM constraint becomes:
\begin{equation}
\mathcal{C}_{rr}(t)\sum |\mathcal{E}_{z,t}(x,y) | /\boldsymbol{\varrho}^{\mathbf{s}_0(t)} /\big( \mathcal{M}_{ag}(\mathbf{s}_0(t))/ \varsigma \big) \leq \lambda(t).
\label{eq:final-fx-certificate}
\end{equation}

\subsection{Integrating FX network feasibility constraints into CEMs}\label{sec:summary}

In summary, solving a CEM with FX network feasibility constraints proceeds in four steps. 
First, an ACOPF solution is computed for a scenario with the same topology, typically over a shorter time horizon, and is used as the reference point.
Second, the certificate~\eqref{eq:final-fx-certificate} is computed based on this reference point.
Third, the resultant linear constraints are embedded within the transport CEM; see Eq.~\eqref{eq_GenX_FXC} below.
Finally, the CEM is solved using these constraints to approximate ACOPF-like investment and dispatch decisions.

Embedding Eq. \eqref{eq_FXC} into the planning problem in Section~\ref{sec_cems}, the GenX model can be reformulated as:
\begin{align} \label{eq_GenX_FXC}
    \min_{x, y} \quad &  \sum_z^\mathcal{Z} \sum_t^T \sum_g^\mathcal{G} (C_\mathcal{I} y+ C_\mathcal{O}x  )   \\
        \text{s. t.} \quad & A_\mathcal{O}x + B_\mathcal{I}y \leq r , \nonumber \\ &F_\mathcal{I}y \leq d,  \nonumber \\ & g_\mathcal{O}(x,y) \leq 0,  
        \nonumber  \\ & \mathcal{C}_{rr}(t)\sum |\mathcal{E}_{z,t}(x,y) | /\boldsymbol{\varrho}^{\mathbf{s}_0(t)} /\big( \mathcal{M}_{ag}(\mathbf{s}_0(t))/ \varsigma \big) \leq \lambda(t), \quad \nonumber  \\ &\forall z \in \mathcal{Z}, \forall t \in \mathcal{T}, \quad x, y \geq 0, x,y \in R.\nonumber 
\end{align}
Note, this FX-CEM does not explicitly model reactive power.
Instead, it aims to constrain active power variables to values that imply admissible reactive power values in practice.

\section{Numerical Experiments and Validation}\label{sec_results}

We evaluate all the proposed formulations on a 
case study of the ISO New England system, featuring three aggregated zones labeled Massachusetts (MA), Connecticut (CT), and Maine (ME).
Hourly active power demand for one year was obtained from public sources and reactive power demand was synthesized (see Section \ref{sec_ReactivePowerGeneration}).
The case study was run with 12 1-week representative periods to allow the ACOPF formulation to solve in  reasonable time. The representative periods were chosen by $k$-means clustering on the demand, solar availability, and wind availability timeseries.
Generation expansion, but not transmission expansion, was allowed. No integer variables were activated. We recognize that the case study does not constitute a fully realistic planning scenario; its purpose is to test how increasingly physics-faithful transmission formulations affect capacity-expansion outcomes in a controlled setting.

All experiments were implemented in Julia \cite{Julia-2017} and executed on a MacBook Pro (2.8\,GHz Intel Core i7, 16\,GB memory). The optimization solver was IPOPT v3.14.19 \cite{wachter2006implementation} with the MA86 linear solver \cite{MA86_report}, MC64 scaling, and an adaptive $\mu$-strategy. 

The case study was solved using five different transmission formulations: (1) the conventional GenX transport formulation; (2) the ACOPF formulation (Section~\ref{sec_ACOPF}; denoted ACOPF in figures and tables); and (3)–(5) the transport formulation with fixed-point–based network feasibility constraints under the tight, medium, and loose bounds (Section \ref{sec:fxconstraints}; denoted by FX-Tight, -Medium, and -Loose). The FX constraints were calculated using a reference point obtained by solving the ACOPF formulation (2) using one representative week;
due to the short time horizon, solving the ACOPF-CEM for one week only added a few seconds to the total runtime.

\begin{figure}[!h]
    \centering  
\includegraphics[width=0.8\linewidth]{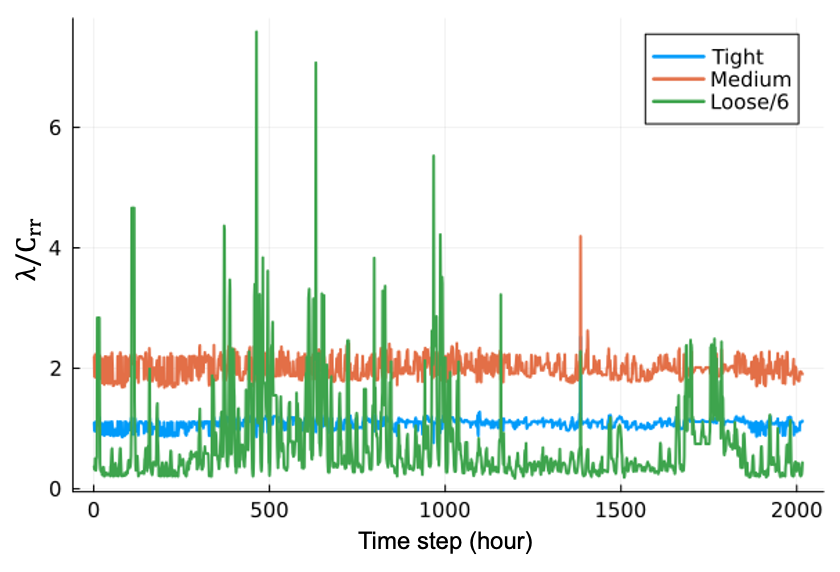 }
\vspace{-1em}
    \caption{ FX feasibility bound ranges for twelve representative weeks of 1 year.}
    \vspace{-1em}
    \label{fig:FXCRanges}
\end{figure}

Fig.~\ref{fig:FXCRanges} presents the computed FX constraints from Eq.~\eqref{eq_FXC}, with tight, medium, and loose bounds. For easier visualization, loose-bound parameters are rescaled to one-sixth of their original values.

\subsection{Generation of AC-Specific Model Data} \label{sec_ReactivePowerGeneration}

Because no year-long hourly reactive power demand data are available for CEMs, we synthesize zonal reactive demand from the active load time series. At each time step, we apply a uniformly sampled power factor between 0.89 and 0.91. We also assign technology-specific reactive power limits: 0.90 power factor for thermal generators, giving symmetric reactive limits of about 50\% of active capacity; 0.98 for storage; and 0.99 for variable renewable generators, reflecting the more limited reactive flexibility of inverter-based technologies \cite{baviskar2020challenges}.

For AC line modeling, we estimate missing parameters—series resistance, reactance, and line charging susceptance—using the standard 380 kV overhead line type (490-AL1/64-ST1A) in the PandaPower library \cite{thurner2018pandapower}, scaled by the line lengths in the three-zone case. We set the system frequency to 60 Hz, base voltages to 345 kV, and zonal voltage magnitudes between 310 and 380 kV. Transformers are included, but their tap ratios are fixed at 1 and angle differences at 0, since transformer effects are not the focus of this paper.

\subsection{Investment Decisions in the Planning Problem} \label{sec:simulation_investments}

As the focus of this study is on  planning outcomes (rather than directly on operations), our metrics of interest are the installed capacities of generators and storage, and the objective cost. We first compare the transport and ACOPF results to evaluate their differences, and then compare ACOPF and FX  
results to see how well the latter approximates the former. ACOPF serves as the ``baseline'' formulation, as it captures AC network physics omitted by the transport formulation.

\subsubsection{Planning with and without ACOPF}

\begin{figure}[!h]
    \centering
    \includegraphics[width=0.8\linewidth]{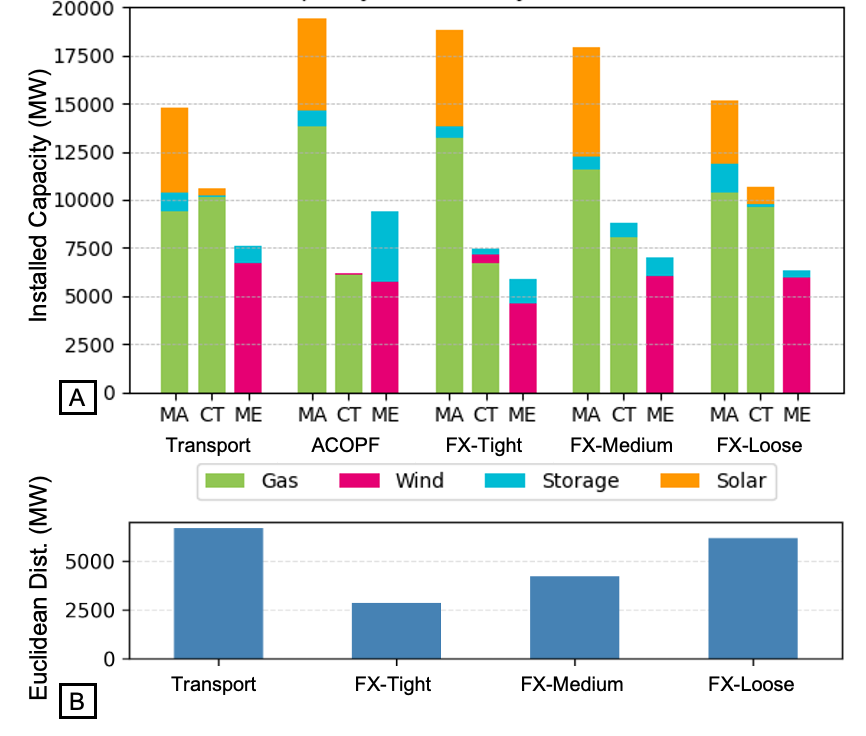}\vspace{-1em}
    \caption{[A] Installed capacities of each generator and storage in the five scenarios. [B] Euclidean distance of the installed capacities in the ACOPF scenario vs. the labeled scenarios.}
\label{fig:InvestmentDecisions}
\end{figure}

Comparing the transport (1) and ACOPF formulations (2) in Fig.~\ref{fig:InvestmentDecisions}.A, we see that the transport formulation significantly underestimates the new capacity required. We quantify this difference by the Euclidean distance between the installed capacities of each generator or storage (Fig.~\ref{fig:InvestmentDecisions}.B), and observe a median absolute difference in individual capacities of 23\% relative to formulation (2). These differences are driven by two factors: additional generators are needed to supply reactive power under formulation (2) (thermal generators are the primary source of reactive power), and the reactive power in the transmission network reduces its ability to transmit active power.

\begin{figure}[!h]
    \centering
    \includegraphics[width=0.85\linewidth]{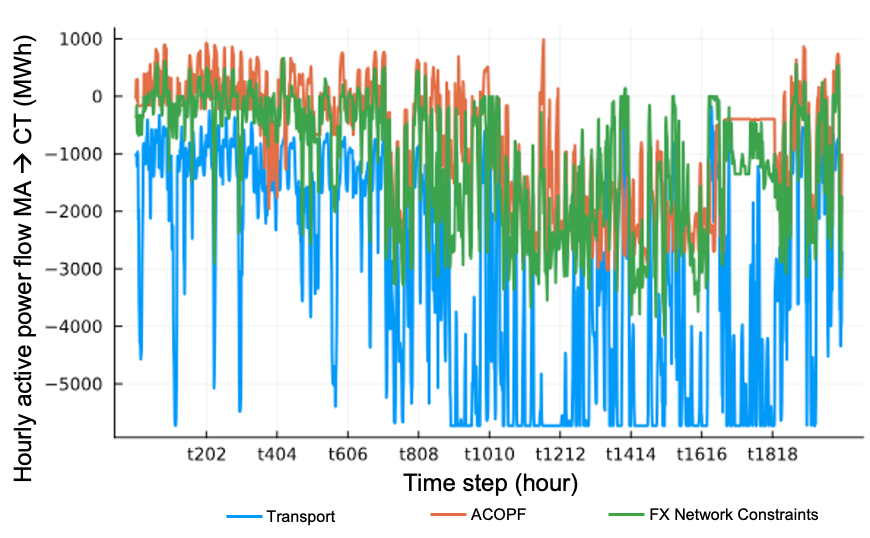}
    \caption{Comparison of active power flow between zone 1 (MA) and zone 2 (CT) under different scenarios.}
    \label{fig:Flows}
\end{figure}

This second point is supported by an analysis of the active power transmission flows. Fig.~\ref{fig:Flows} reports the active power flow between zone 1 (MA) and zone 2 (CT) under formulations (1)-(3). The capacity limit of the transmission line is 5730\,MW. The black curve in Fig.~\ref{fig:Flows} shows that the network frequently hits the transmission limit under formulation (1). Under formulation (2), the maximum of the active power flow on the line is 3738\,MW due to the reactive power flow on the line, while the maximum magnitude of the apparent power is $\approx$~5277\,MW.
The more detailed treatment of reactive power generation and transmission reduces the grid's ability to transmit active power, requiring more local generation.

\subsubsection{Planning with FX constraints}

\begin{figure}[!h]
    \centering
\includegraphics[width=0.8\linewidth]{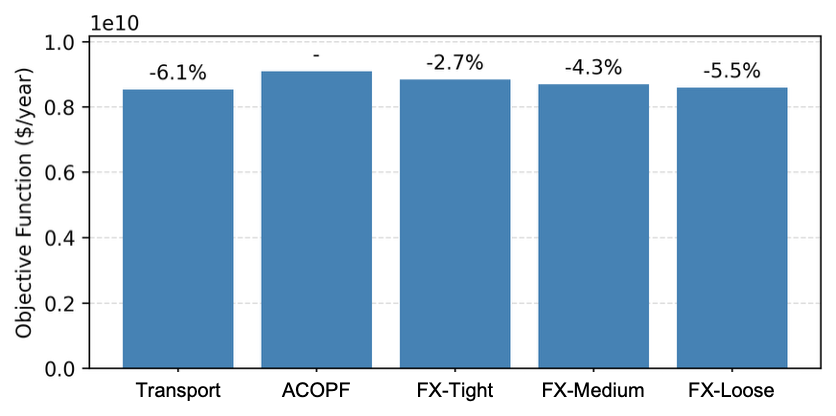}\vspace{-1em}
    \caption{Comparison of total costs for each case. The \% difference from the ACOPF result is given above each case}
    \label{fig:costs}
\end{figure}

Fig.~\ref{fig:InvestmentDecisions} compares formulations (3)–(5) with formulations (1)–(2) using Euclidean distance. We can see that the FX constraints give planning results closer to those of the ACOPF-CEM than the transport formulation. Formulation (3), with the tightest bounds, yields the closest decisions, with 10\% median absolute difference in individual generator capacities. The looser bounds have greater errors:~16\% and 30\% for the medium and loose bounds, respectively. In this respect, the loose FX bound is worse than the transport formulation. Total system costs are similar across formulations, as shown in Fig.~\ref{fig:costs}, with ACOPF highest, transport lowest, and the FX cases ordered by bound tightness.

Overall, these results support the effectiveness of the proposed FX approach. Even though our FX-CEM does not include reactive power explicitly, the FX network feasibility constraints capture the effects of AC physics on network congestion. This is evident in Fig.~\ref{fig:Flows}; the active power line flows under formulation (3) are similar to those under formulation (2). Explicitly representing reactive power in the FX-CEM may improve its accuracy further.

\subsection{Comparison of Model Dimensions \& Computational Costs}
Table \ref{table:ModelSizeTime_compare} compares the model size and computational effort required across the three transmission formulations. As expected, the transport formulation is the smallest and fastest given that it is linear (0 non-linear expressions; 0 Hessian nonzeros). Adding ACOPF increases the number of variables and constraints by 50\%. This includes 16K non-linear expressions and dense Jacobians/Hessian, increasing solve time by 4$\times$. The transport formulation with FX constraints does not introduce new variables and requires only a small number of new constraints (3\% of the original model), increasing  runtime by 5--20\%. This study's version of GenX with FX constraints includes non-linear expressions due to the use of absolute values. These have been linearized without loss of generality in subsequent versions, making the model a MILP. 

\begin{table}[!h]
\centering
\caption{Comparison of Transmission Formulations}
\scalebox{0.8}{
\begin{tabular}{l|l|l|l}
\hline\hline
Formulation & Transport & \begin{tabular}[c]{@{}l@{}} ACOPF\end{tabular}                   & \begin{tabular}[c]{@{}l@{}}Transport \\ \& FX\end{tabular}                                    \\ \hline\hline
Total time (seconds) & 154.08 - 179.45  & 460.06 - 643.53  &  182.02- 191.72   \\
\hline
IPOPT time (seconds) &150.32-170.33 & 443.79 - 627.89  &  140.80- 176.34   \\
\hline
Total iterations &198-227 & 353-422 &  211 -253  \\
\hline
\begin{tabular}[c]{@{}l@{}}Total number  \\ of variables \end{tabular}  & 100813                                                                                           & 153229                                                                         & 100813 \\ \hline
\begin{tabular}[c]{@{}l@{}}Total number \\ of constraints\end{tabular}            & 209683                                           & 292339                                                                    & 215731                           \\ \hline
\begin{tabular}[c]{@{}l@{}}Total number \\ of non-linear \\ expressions\end{tabular}      & 0 & \begin{tabular}[c]{@{}l@{}}16128 \end{tabular}                  & 6048                                  \\ \hline
\begin{tabular}[c]{@{}l@{}}Number of nonzeros \\ in equality\\ constraint Jacobian\end{tabular}              & 94752             & \begin{tabular}[c]{@{}l@{}} 225792\end{tabular} & 94752    
\\ \hline
\begin{tabular}[c]{@{}l@{}} Number of nonzeros \\ in inequality \\constraint Jacobian\end{tabular}              & \begin{tabular}[c]{@{}l@{}} 185446  \end{tabular}             & \begin{tabular}[c]{@{}l@{}}288313 \end{tabular} & 193510
\\ \hline
\begin{tabular}[c]{@{}l@{}} Number of nonzeros \\ in Lagrangian  \\Hessian\end{tabular}              & \begin{tabular}[c]{@{}l@{}} 0  \end{tabular}             & 177408  & 10080
\\ \hline
\end{tabular}}
\label{table:ModelSizeTime_compare}
\vskip -1em
\end{table}

\subsection{Robustness of Computed Parameters}
The results of the first case study with formulations (3)–(5) relied on fixed points obtained from the same case study under formulation (2). In this section, we test the robustness of the FX approach by running a second case study in which formulations (3)–(5) use fixed points calculated from different, third and fourth case studies under formulation (2).

\begin{figure}[!h]
    \centering
    \includegraphics[width=0.8\linewidth]{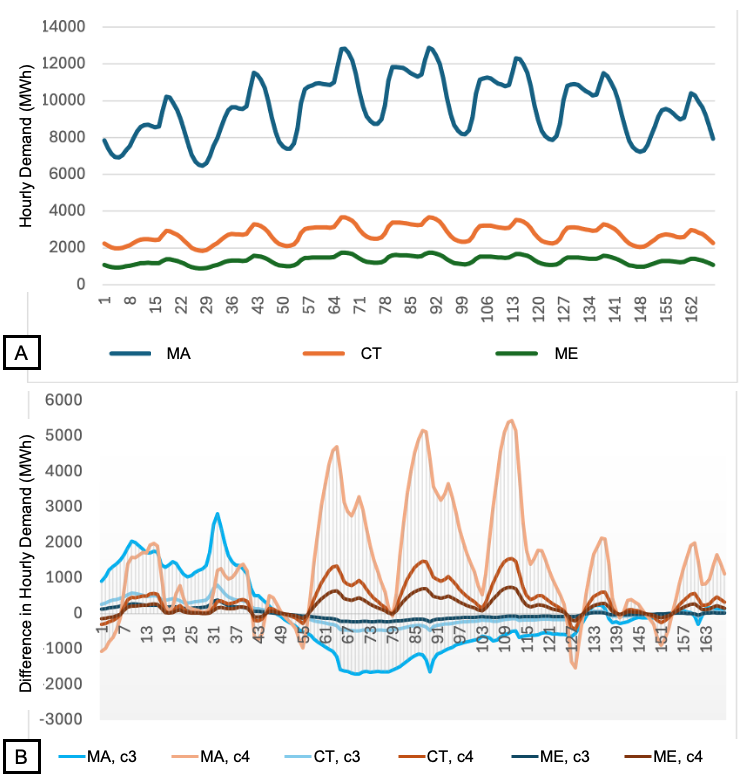}
    \caption{[A] Active power demand of each zone in the ISO-NE case study for the representative week used for as the fixed point in formulation (3). [B] Difference in zone demand for the third and fourth case studies}
    \label{fig:Demand}
\end{figure}

The second case study is the same as the first but with lower fuel prices. The third and fourth case studies have the same lowered fuel prices and also altered demand timeseries. The original representative demand timeseries in the first and second case studies is shown in Fig.~\ref{fig:Demand}.A. Figure \ref{fig:Demand}.B shows the differences in the zonal demand for the third and fourth case studies, denoted $c$ and $c$ respectively. The third case study has lower peak demand during the peak hours of the first case study, with similar or modestly lower levels elsewhere. The fourth case study has uniformly higher load.

\begin{figure}[!h]
    \centering
    \includegraphics[width=0.8\linewidth]{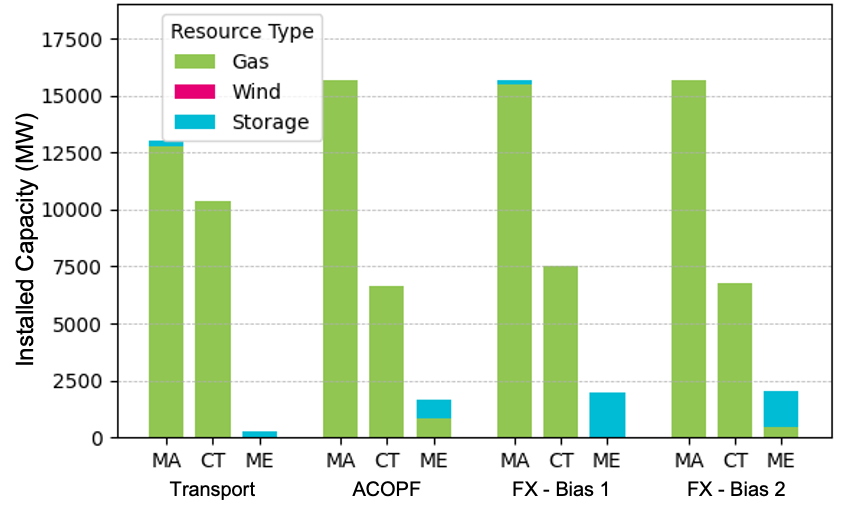}
    \caption{Installed capacities for the second ISO-NE case study with low fuel costs. The two sets of FX network constraints were calculated from the third and fourth case studies, with the demand alterations shown in Figure \ref{fig:Demand}.B.}
    \vspace{-1em}
    \label{fig:Robustness}
\end{figure}

We solve the second case study using four transmission formulations: the original transport formulation GenX, GenX with ACOPF constraints, and two versions of GenX with FX network feasibility constraints, using the tight bounds and reference points taken from the third and fourth case studies solved using formulation (2).

Fig.~\ref{fig:Robustness} compares investment decisions across the four formulations. Both FX formulations lead to investment results closer to the ACOPF reference than the transport formulation, indicating that parameters calibrated on one dataset can generalize to out-of-sample demand profiles. This robustness is valuable for planning under uncertainty and for very long time periods.

\section{Conclusion}

This paper demonstrates that incorporating AC power network physics changes optimal capacity expansion planning decisions. We compared a capacity expansion model using a transport formulation of transmission with one incorporating ACOPF. 
In our test case, we demonstrated that the latter installed more new capacity and shifted investments toward local resources, due to reduced active power transfer capability when the modeled network must also transmit reactive power. However, the ACOPF-based formulation took four times longer to run than the transport formulation. To overcome this scaling challenge, we introduced fixed-point–based network feasibility constraints as an effective surrogate that embeds ACOPF effects within a linear planning model. 

The proposed FX constraints achieve a favorable accuracy–cost trade-off: they closely reproduce ACOPF investment patterns with near-transport dimensionality and far lower runtimes than ACOPF. This study shows that effective FX constraints can be calculated using only a few representative weeks from a case study, allowing them to be used in planning models with long time horizons. We have also shown that FX parameters calculated on one set of demand data can be used for other demand data on the same network.
\textcolor{black}{Future work includes exploring the implications of different choices of nominal and reference points, incorporating reactive power explicitly within the FX formulation, handling of transmission expansion, a continuous choice of FX bounds, further improving FX constraints' robustness under different scenarios, and conducting analysis on additional power systems.}



\bibliographystyle{IEEEtran}
\bibliography{ieeebib}
\end{document}